\documentclass[12pt]{article}
\usepackage{amsmath,amssymb, amsfonts, amsthm,amscd}
\usepackage[T2A]{fontenc}
\usepackage[cp1251]{inputenc}
\usepackage[russian,ukrainian,english]{babel}
\usepackage{graphicx}

\newtheorem*{proposition*}{Proposition}
\newtheorem*{lemma*}{Lemma}
\newtheorem*{theorem*}{Theorem}
\newtheorem*{corollary*}{Corollary}

\newtheorem*{definition*}{Definition}
\newtheorem*{remark*}{Remark}
\newtheorem*{example*}{Example}

\begin{document}

\begin{center}
{\LARGE Diagonal reduction of matrices over commutative semihereditary  Bezout rings }
\end{center}
\vskip 0.1cm \centerline{\large \large B.~Zabavsky,  A.~Gatalevych}
\vskip 0.5cm

\centerline{Ivan Franko National University of Lviv,
 1 Universytetska str.,Lviv, 79000, Ukraine}

\centerline{gatalevych@ukr.net}
\medskip

Abstract. It is proven  that every commutative semihereditary Bezout ring in which any
regular element is Gelfand (adequate), is an elementary divisor ring.
\vskip 0.5cm
MR(2010)Subject Classification 06F20, 13F99

\vskip 0.5cm

All rings considered will be commutative and have identity. Recently there has been some interest
in the polynomial ring $R[x]$, where $R$ is a von Newman regular ring. Such a ring is a Bezout ring [6],
semihereditary ring [6], and so Hermite ring [4]. Thus, it is natural to ask whether or not $R[x]$
is an elementary divisor ring. This question is answered affirmative in [12]. It is an open problem
whether or not every Bezout domain is an elementary divisor ring and more generally: whether or not every semihereditary Bezout ring is an elementary divisor ring.

In a recent paper by Shores [12], we obtain a complete characterization of  semihereditary  elementary divisor ring through its homomorphic
images. This result will be pivotal in our work.

{\tt Theorem 1.} (see [12]) ("Shores" test) Let $R$ be a semihereditary Bezout ring. Then  $R$ is an elementary divisor ring iff $R/aR$ is an elementary divisor ring for all nonzero divisor (regular) elements $a\in R$.

We introduce the necessary definitions and facts.

By a Bezout ring we mean a ring in which all finitely generated ideals are principal.
 An $n$ by $m$ matrix $A=(a_{ij})$ is said to be diagonal if
 $a_{ij}=0$ for all $i \neq j$. We say that a matrix $A$ of the
 dimension $n$ by $m$ admits a diagonal reduction if there exist
 invertible matrices $ P \in GL_n(R),Q \in GL_m(R) $ such that $ PAQ $
 is a diagonal matrix. We  say that two matrices $A$ and $B$
 over a ring $R$ are equivalent  if there
 exist invertible matrices $P,Q$ such that $B=PAQ$.
 Following Kaplansky [4], we say that if every matrix over $R$ is
 equivalent to a diagonal matrix $(d_{ii})$ with the property that
 every $(d_{ii})$ is a  divisor of $ d_{i+1,i+1}$, then $R$ is an
 elementary divisor ring.
A ring $R$ is to be a  Hermite ring if every
 $1\times 2 $ matrix over $R$ admits diagonal reduction.

A row $(a_1; a_2;\dots ; a_n)$ over a ring R is called unimodular if $a_1R + a_2R\dots + a_nR = R$. If
 $(a_1; a_2;\dots ; a_n)$ is a unimodular $n$-row over a ring $R$, then we say that $(a_1; a_2;\dots ; a_n)$ is reducible
if there exists $(n-1)$-row $(b_1; b_2;\dots ; b_{n-1})$ such that the $(n-1)$-row
$(a_1+a_nb_1; a_2+a_nb_2;\dots ; a_{n-1}+a_nb_{n-1})$ is unimodular. A ring $R$ is said to have a stable range $n$ if $n$ is the least positive integer
such that every unimodular $(n+1)$-row is reducible.

Obviously, an elementary divisor ring is a Hermite ring, and it is easy to see that a a Hermite ring is a Bezout ring. Examples that neither
implication is reversible are provided by Gillman and Henriksen [2].  In [13], it was proved that a commutative Bezout ring is Hermite iff it is a ring of stable range 2.

A von Neumann regular ring is a ring $R$ such that for every $a\in R$ there exists an $x \in R$ such that $a=axa$.
A  ring $R$ is  said  to  be  semiregular  if $R/J(R)$ is  regular  and idempotents  lift  modulo $J(R)$,
where $J(R)$ denotes  the  Jacobson  radical  of $R$ [9].

An element $a\in R$ is called clean if $a$ can be written as the sum of  a unit and an idempotent. If
each element of $R$ is clean, then we say $R$ is a clean ring [9].

A clean ring is a Gelfand ring. Recall that a ring is called a Gelfand ring if whenever
$a+b=1$ there are $r, s \in R$ such that $(1+ar)(1+bs)=0$ [7].

 A ring is called a $PM$-ring if
every prime ideal is contained in a unique maximal ideal [1].

It had been asserted that a commutative ring is a Gelfand ring iff it is a $PM$-ring [1].

A ring is called a $PM^{\ast}$-ring if
every nonzero prime ideal is contained in a unique maximal ideal [14].

An element $a$ of the ring $R$  is called an atom if from decomposition $a=b\cdot c$ follows that $b$ or $c$ is a unit of $R$.

The group of units of a ring $R$ will be denoted by $U(R)$, the Jacobson radical of a ring $R$ will be denoted by $J(R)$.
For a ring $R$, $spec(R)$ ($mspec(R)$) denotes the collection of prime (maximal) ideals of $R$, and $spec(a)=\{ P\in spec(R)| a\in P\}$,
$mspec(a)=\{ M\in mspec(R)| a\in M\}$.

We start with trivial statements that are of a technical nature.

{\tt Proposition 2.} Let $R$ be a Bezout ring. Then element $a\in R$ is an atom iff the factor ring $R/aR$  is a field.

{\it Proof.} Let $a\in R$ be an atom. Denote $\bar b=b+aR$ for some $b\in R$ and $\bar R=R/aR$. Let $\bar b\neq\bar 0$, then $b\notin aR$.
Since $R$ is a Bezout ring, we  obtain $aR+bR=dR$ and $a=a_0d, b=b_0d$, $au+bv=d$ for some $d, a_0, b_0, u, v \in R$. Since $a$ is an atom, we have $d\in U(R)$ or $a_0\in U(R)$. If  $d\in U(R)$, we have $\bar d\in U(\bar R)$ and then $\bar b\in U(\bar R)$, where $\bar d=d+aR$. If $a_0\in U(R)$, we have $b=aa_0^{-1}b_0\in aR$ and this is  impossible, since $\bar b\neq \bar 0$.
So we proved that $\bar R$ is a field.

If  $\bar R$ is a field, then $aR\in mspec(R)$ and, obviosly, $a$ is an atom in $R$. Proposition is proved.

Now let $R$ be a  ring and let $a\in R\setminus{0}$. Let us denote $\bar R=R/aR$ and $\bar b=b+aR$.

{\tt Proposition 3.} Let $R$ be a ring. Then $\bar b\in U(\bar R)$ iff $aR+bR=R$.

 {\it Proof.} Let $\bar b\in U(\bar R)$, then $\bar b\bar v=\bar 1$, where $\bar v=v+aR$, and $bv-1\in aR$ i.e. $aR+bR=R$.

 If $aR+bR=R$, then $au+bv=1$ for some $u, v \in R$ and $\bar b\bar v=\bar 1$.

Proposition is proved.

Regarding the concept of a stable range of the ring, one can distinguish a stable range of 1 and 2.

A ring $R$ is said to have stable range 1 if for all $a,b\in R$
such that $aR+bR=R$, there exists $y\in R$ such that $(a+by)R=R.$

 A ring $R$ is said to have stable range 2 if for all $a,b, c\in R$
such that $aR+bR+cR=R$, there exists $x, y \in R$ such that $(a+cx)R+(b+cy)R=R$.

 {\tt Proposition 4.}
Let $R$ be a ring of stable range 2 and let $aR+bR=dR$. Then there exist  elements $a_0, b_0\in R$ such that
$a=a_0d, b=b_0d$ and $a_0R+b_0R=R$.

{\it Proof.} Since $aR+bR=R$, then $a=a_1d, b=b_1d$, $au+bv=d$ for some $d, a_1, b_1, u, v \in R$. Then $d(1-a_1u-b_1v)=0$.
Let us denote $c=1-a_1u-b_1v$, then $dc=0$ and $a_1R+b_1R+cR=R$. Since $R$ is a ring of stable range 2 we have
$(a_1+cx)R+(b_1+cy)R=R$ for some elements $x, y \in R$.
Let us denote $a_0=a_1+cx, b_0=b_1+cy$. It is obvious that $a_0d=(a_1+cx)d=a_1d+cxd=a_1d=a$, $b_0d=b$ and $a_0R+b_0R=R$.

Proposition is proved.

{\tt Proposition 5.} Let $R$ be such a ring that for any elements $a, b \in R$ there exist such elements
$d, a_0, b_0 \in R$ that $a=da_0, b=db_0$ and $a_0R+b_0R=R$. Then $R$ is a Bezout ring of stable range 2.

{\it Proof.} Since $a_0R+b_0R=R$ and $a=da_0, b=db_0$ we have $a_0u+b_ov=1$ and $au+bv=d$. It implies that $d\in aR+bR$, i.e. $dR\subset aR+bR$.
Since $a=da_0, b=db_0$, we have $aR+bR \subset dR$, and therefore
 $aR+bR=dR$, i.e. $R$ is a Bezout ring.

Next we show that the ring $R$ is a ring of stable range 2. Let $aR+bR+cR=R$. According to the restrictions imposed on $R$, we have
$a=da_0, b=db_0$ and $a_0R+b_0R=R$ for some elements $a_0, b_0\in R$. Since  $a_0R+b_0R=R$ we have  $a_0u+b_0v=1$ for some elements $u, v\in R$.
Then $(a-cv)(-b_0)+(b+cu)a_0= -ab_0+ba_0+c(a_0u+b_0v)$. Since $ab_0=ba_0$ and $a_0u+b_0v=1$, we have
$(a-cv)(-b_0)+(b+cu)a_0=c$. In addition,  we obtain $(a-cv)u+(b+cu)v=au-cvu+bv+cuv=au+bv=d$. Since $aR+bR+cR=R$ and $aR+bR=dR$ we obtain
$dR+cR=R$. Then $(a+c(-v))R+(b+cu)R=R$ i.e.,  $R$ is a ring of stable range 2.

Proposition is proved.

Mc Adam S. and  Swan R. G. [5] studied comaximal factorization in commutative rings. Following them, we give the following definitions.

{\tt Definition 1.} A nonzero element $a$ of a ring $R$ is called inpseudo-irreducible if for any representation $a=b\cdot c$ we have $bR+cR=R$.

{\tt Definition 2.} An element $a$ of a ring $R$ is called pseudo-irreducible if for any representation $a=b\cdot c$, where $b, c \notin U(R)$, we have $bR+cR\neq R$.

For example, for the ring of integers $\mathbb{Z}$: 6 is inpseudo-irreducible element, and 4 is pseudo-irreducible element.

{\tt Definition 3.} A nonzero element $a\in R$ is called a regular element if $a$ is not  a zero divisor of $R$.

{\tt Proposition 6.} Let $R$ be a Bezout ring and $a$ is a regular inpseudo-irreducible element of $R$.
Then $R/aR$ is a reduced ring.

{\it Proof.}  Let $\bar R=R/aR$ and $\bar x=x+aR, \bar x^n=\bar 0$ for some  $n\in \mathbb{N}$. Obviously that $\bar x\notin U(\bar R)$.
Let $xR+aR=dR$. Since $\bar x\notin U(\bar R)$ by Proposition 3,  we have that $d\notin U(R)$. In this case $x=dx_0, a=da_0$ and $xu+av=d$ for some elements $x_0, a_0, u, v \in R$. Then $d(1-x_0u-a_0v)=0$. Since $d\neq 0$ and $d$ is a divisor  of the regular element $a$ we have $x_0u+a_0v=1$. Since $a$ is an inpseudo-irreducible element and $a=da_0$, we obtain that $a_0R+dR=R$. As we have already proved that $a_0R+x_0R=R$ and $a_0R+dR=R$ we  obtain,  that $a_0R+xR=R$. Whereas $\bar x^n=\bar 0$, we have $x^n=at$ for some element $t\in R$.
By $x=dx_0, a=da_0$ and $x^n=at$ we have that $d^nx_0^n=da_0t$, and  $d(d^{n-1}x_0^n-a_0t)=0$. Since $d\neq 0$ and $d$ is a divisor  of the regular element $a$, we have $d^{n-1}x_0^n-a_0t=0$.

From the condition   $xR+a_0R=R$, it follows that  $x^{n-1}R+a_0R=R$ and $x_0R+a_0R=R$. We have $d^{n-1}x_0^n=a_0t$, i.e.  $x^{n-1}x_0=a_0t$.
Since $x^{n-1}x_0R+a_0R=R$, we have $a_0\in U(R)$, i.e. $xR\subset aR$. It follows that $R/aR$ is a reduced ring.

Proposition is proved.

{\tt Proposition 7.} Let $R$ be a Bezout ring and $a$ is a regular  element of $R$. Then an annihilator of each element of
 $R/aR$ is a principal ideal.

{\it Proof.} Let $\bar R=R/aR$ and let $\bar b=b+aR$ be any nonzero element of $\bar R$. Let $Ann \bar b=\{\bar x\in \bar R| \bar x \bar b=\bar 0\}$, where
$\bar x=x+aR$.  As $\bar x\neq \bar 0$ it is clear that $bx=ay$ for some element $y\in R$. Since $R$ is a Bezout ring, $aR+bR=dR$ and $a=a_0d, b=b_0d$, $au+bv=d$
for some elements $a_0, b_0, u, v \in R$.   Then $d(1-a_0u-b_0v)=0$. Whereas $d$ is a nonzero divisor of the regular element $a$, we have $a_0u+b_0v=1$.
Then we get the following relationships $b_0x=a_0y$, $a_0ux+b_0vx=x$ and $x=a_0ux+a_0yv$.  Then $x\in a_0R$. From $a_0b=ab_0$ it follows that $\bar a_0\in Ann\bar b$.
From the arbitrariness of the element $\bar x$ and from the conditions $x\in a_0R$ and $\bar a_0\in Ann  \bar b$, we obtain that $Ann\bar b=\bar a_0\bar R$.
Moreover, $Ann\bar b=\bar a_0\bar R$.

Proposition is proved.

{\it Remark 1.} Let $R$ be a Bezout ring and $a\in R$ is a regular  element. Let $\bar R=R/aR.$  By Proposition 7 we have that $Ann\bar b=a_0\bar R$ for any element $\bar b=b+aR$, where $aR+bR=dR$, and $a=a_0d, b=b_0d$, $au+bv=d$. Then obviously $Ann\bar a_0=\bar d\bar R$ and since
$au+bv=d$, we obtain that $Ann(Ann\bar b)=\bar b\bar R$ for any element $\bar b\in \bar R$.

{\tt Proposition 8.} Let $R$ be a Bezout ring and $a$ is a regular element of $R$.
Then $R/aR$ is a reduced ring iff $a$ is an inpseudo-irreducible element.

{\it Proof.}  Let $\bar R=R/aR$  be  a reduced ring.
Let $a=b\cdot c$ and $bR+cR=dR$, where $b=b_ 0d, c=c_0d$ for some elements $b_0, c_0 \in R$. Consider

$$(bc_0)^2=bcc_0b_0=ac_0b_0\in aR.$$

Since $\bar R$ is reduced, we have $bc_0=at$ for some element $t\in R$. It follows that $bc=a=adt$ and by regularity of the element $a$ we obtain $dt=1$, i.e. $a$ is an inpseudo-irreducible element.

By Proposition 6, Proposition is proved.

From the Propositions 6, 7, 8  we will obtain the following results.

{\tt Theorem 9.} Let $R$ be a Bezout ring of stable range 2.
A regular element $a\in R$ is inpseudo-irreducible iff  $R/aR$ is a von Neuman regular ring.

{\it Proof.} Let us prove the necessity, namely, we will show that $R/aR$ is a von Neuman regular ring.

By Proposition 6, we have that $\bar R=R/aR$ is a reduced ring. By Proposition 7, we have that in $\bar R$
annihilator $Ann\bar b$ of each element $\bar b \in \bar R$ is a principal ideal, i.e. $Ann\bar b=\bar{\alpha}\bar R$.
Let us see that $\bar b\bar R\cap\bar{\alpha}\bar R=\{\bar 0\}$. Indeed, let $\bar k\in\bar b\bar R\cap\bar{\alpha}\bar R$,
i.e. $\bar k=\bar b\bar t=\bar{\alpha}\bar s$ for some elements $\bar t, \bar s\in \bar R$. Then
$\bar k^2=\bar b\bar t\bar{\alpha}\bar s=\bar b\bar{\alpha}\bar t\bar s=\bar 0\bar t\bar s=\bar 0$.
Since $\bar R$ is reduced, we obtain $\bar k=\bar 0$, so $\bar b\bar R\cap\bar{\alpha}\bar R=\{\bar 0\}$.

Next, note that $\bar R$ is a commutative Bezout ring of stable range 2. Let $\bar b\bar R+\bar{\alpha}\bar R=\bar{\delta}\bar R$.
Then by Proposition 4, we have $\bar b=\bar b_0\bar{\delta}, \bar{\alpha}=\bar{\alpha}_0\bar{\delta}$ and $\bar b_0\bar R+\bar{\alpha}_0\bar R= \bar R$ for some elements $\bar b_0, \bar{\alpha}_0\in \bar R$.
Since $\bar b_0\bar{\alpha}\bar b\bar{\alpha}_0\in \bar b\bar R\cap\bar{\alpha}\bar R$, then $\bar b_0\bar{\alpha}=\bar b\bar{\alpha}_0=\bar 0$.
Since $\bar b\bar{\alpha}_0=\bar 0$, then $\bar{\alpha}_0\in Ann\bar b=\bar{\alpha}\bar R$,
i.e. $\bar{\alpha}_0=\bar{\alpha}\bar t$ for some element $\bar t\in \bar R$. Then
$$\bar b_0\bar{\alpha}_0=\bar b_0\bar{\alpha}\bar t=\bar b_0\bar{\delta}\bar{\alpha}_0\bar t=\bar b\bar{\alpha}_0\bar t\in
\bar b\bar R\cap\bar{\alpha}\bar R=\{\bar 0\}.$$

By Remark 1, we have $Ann\alpha = \bar b\bar R$. Let us consider $\bar b_0\bar{\alpha}=\bar 0$, then $\bar b_0\in Ann\bar{\alpha}=\bar b\bar R$, i.e. $\bar b_0\in \bar b\bar R$. Thus we proved
$\bar{\alpha}_0\bar R\subset \bar{\alpha}\bar R$ and $\bar b_0\bar R\subset \bar b\bar R$. Therefore

$$\bar R=\bar{\alpha}_0\bar R+\bar b_0\bar R\subset \bar{\alpha}\bar R+\bar b\bar R \Rightarrow \bar{\alpha}\bar R+\bar b\bar R=\bar R.$$

Then $\bar b\bar u+\bar{\alpha}\bar v=\bar 1$ for some elements $\bar u, \bar v\in \bar R$. Since $\bar b\bar{\alpha}=\bar 0$, we obtain
$\bar b^2\bar u=\bar b$, i.e. $\bar b$ is a von Neumann regular element and $\bar R$ is a von Neumann regular ring.

Since the commutative von Neumann regular ring is reduced, sufficiency follows from Proposition 8.

{\tt Definition 4.} An element $a\in R$ is called an element of almost stable range 1 if the ring $R/aR$  is a ring of stable range 1.

Since the commutative von Neumann regular ring is a ring of stable range 1 [16], we obtain the following result.

{\tt Corollary 10.} Let $R$ be a Bezout ring of stable range 2.
A regular inpseudo-irreducible element is an element of almost stable range 1.

Now consider the case of an pseudo-irreducible element.

As it is well-known, a ring is said to be indecomposable if it cannot
be decomposed into a direct sum of two or more non-trivial ideals.
A ring $R$ is an indecomposable ring if and only if 1 is the only non-zero idempotent of $R$.

{\tt Theorem 11.} Let $R$ be a Bezout ring of stable range 2.
A regular element $a\in R$ is pseudo-irreducible iff  $R/aR$ is an indecomposable ring.

{\it Proof.} Let $a\in R$ be a regular pseudo-irreducible element, and let $\bar e=e+aR$, $\bar e^2=\bar e$.  Then $\bar e(\bar 1-\bar e)=\bar 0$, i.e.
$e(1-e)=at$ for some element $t\in R$. Assume that $\bar e\neq \bar 0$, i.e. $e\notin aR$, and $\bar e\neq \bar 1$, i.e. $eR+aR\neq R$ and $e\notin aR$. Let $eR+aR=dR$. By proposition 4, we have $e=e_0d$, $a=a_0d$, and $e_0R+a_9R=R$, where $d\notin U(R)$ and $a_0\notin U(R)$. Since $a$ is a regular element, we have $e_0(1-e)=a_0t$. From the condition  $e_0R+a_0R=R$, it  follows that $e_0u+a_0v=1$ for some elements $u, v \in R$. Then
$(1-e)e_0u+(1-e)a_0v=1-e$ and by $e_0(1-e)=a_0t$ we have $a_0(ut+(1-e)v=1-e$, i.e. $a_0\alpha +e=1$ for some element $\alpha \in R$.
Then $a_0R+eR=R$ and $a_0R+dR=R$, i.e. we obtain that $a=da_0$ where $d\notin U(R), a_0\notin U(R)$ and $dR+a_0R=R$.  We obtain  the contradiction, therefore $\bar R$ is an indecomposable ring.

Let $R/aR$ be an indecomposable ring. Suppose that $a$ is not an pseudo-irreducible element, i.e. $a=bc$, where $b\notin U(R), c\notin U(R)$ and
$bR+cR=R$ for some elements $b, c \in R$, and
$bu+cv=1$ for some elements $u, v \in R$. Then $\bar b^2\bar u=\bar b$, where $\bar b=b+aR, \bar u=u+aR$. Since $\bar R$ is an indecomposable ring and by $\bar b\bar u=\bar e$, $\bar e^2=\bar e$, we have $\bar b\bar u=\bar 0$ or $\bar b\bar u=\bar 1$.
If $\bar u\bar u=\bar 0$  then $bu=ct$ for some element $t\in R$. It follows  that $at+cu=1$, i.e. $c\in U(R)$ but this is impossible.
If $\bar b\bar u=\bar 1$  then by Proposition 3, it follows that $aR+bR=R$. Since $a=bc$, we obtain that $b\in U(R)$ and this is impossible.
Thus we proved that $a$ is an pseudo-irreducible element.

Theorem is proved.

A basic property of clean rings is that any homomorphic image
of a clean ring is again clean. This leads to the definition of a neat ring.

{\tt Definition 4.} A ring $R$ is
a neat ring if every nontrivial homomorphic image of $R$ is clean [7].

Let $R$ be a ring in which for each regular element $a\in R$, the ring $R/aR$ is a clean ring (for example $R$ is a neat Bezout ring). Then by Theorem 11, we have that the  regular element $a$ is pseudo-irreducible iff $R/aR$ is a local ring, i.e. the element $a$  is contained in a unique maximal ideal. For a principal ideal domain $R$, the pseudo-irreducible elements will be elements of the form $p^n$, where $p$ is an atom of $R$ and $n\in \mathbb{N}$.
But this will not always be the case. For example, let

$$ R=\{z_0+a_1x+a_2x^2+ \dots +a_nx^n+\dots | z_0\in \mathbb{Z}, a_i\in \mathbb{Q}\}.$$

Here the  element $x$  is pseudo-irreducible, but the ring $R/xR$ is not local.

{\tt Definition 5.} Let $R$ be a Bezout ring. An element $a\in R$  is called
an adequate element  if   for every $b\in R$  there exist such $r, s\in R$  that

(i) $a = rs$,

(ii) $rR+bR=R$, and

(iii) for each non-unit divisor $s'$ of $s$, we have $s'R+bR\neq R$.

If every non-zero element of $R$ is adequate, the ring $R$ is called an adequate ring.

{\tt Theorem 12.} Let $R$ be a Bezout ring of stable range 2.
A regular element $a\in R$ is an adequate  element iff  $R/aR$ is a semiregular ring.

{\it Proof.} Let  $\bar b=b+aR$ be a non-zero and non-invertible element of $\bar R=R/aR$.
Since $a$ is an  adequate  element of $R$, we have  $a = rs$, where $rR+bR=R$, and $s'R+bR\neq R$
 for every non-unit divisor $s'$ of $s$. Note that $rR+sR=R$. Indeed,  otherwise $rR+sR=hR\neq R$, then, according to
 the definition of an adequate element $a$, we, on the one hand, obtain that $hR+bR=R$ ($h$ is a divisor of $r$),  and on the other hand,
 $hR+bR\neq R$ ($h$ is non-unit divisor of element $s$), which is imposible.  This means that $rR+sR=R$ and $ru+sv=1$ for some elements $u, v \in R$. Then $\bar r\bar R+\bar s\bar R=\bar R$, where $\bar r=r+aR$, $\bar s=s+aR$. By analogy, we get that
 $\bar r\bar R+\bar b\bar R=\bar R$.

 Denote $\bar e=\bar s\bar v$, then $\bar e^2=\bar s\bar v(\bar 1-\bar r\bar u)=\bar s\bar v=\bar e$ and then
 $\bar 1-\bar e=\bar r\bar u$.
 Since $\bar r\bar R+\bar b\bar R=\bar R$ and $\bar r^2\bar u=\bar u$, we have $(\bar 1 - \bar e)\bar R+\bar b\bar R=\bar R$. It follows that
 $\bar e\in \bar b\bar R$ and $mspec(\bar b)\subset mspec(\bar e)$.  Let there exists $\bar M \in mspec(\bar 1- \bar e)$ and $\bar b \notin \bar M$. Then $\bar M+\bar b\bar R=\bar R$, i.e. $\bar m+\bar b\bar t=\bar 1$, where $\bar m\in\bar M, \bar t\in \bar R$.

Let $\bar e\bar R+\bar m\bar R=\bar d\bar R$. From $\bar e, \bar m\in \bar M$ it follows that $\bar d\in\bar M$.
Since $\bar d$ is a non-unit divisor of $\bar s\bar v$ and from $\bar s\bar v=\bar e$ and $\bar s^2\bar v=\bar s$, it follows that $\bar d$
is a  non-unit divisor of $\bar s$. Then $\bar b\bar R+\bar d\bar R\neq \bar R$. But
$\bar R=\bar b\bar R+\bar m\bar R\subset\bar b\bar R+\bar d\bar R\neq \bar R$, that is impossible. Hence,  $mspec(\bar b)=mspec(\bar e)$.

We will prove that $\bar b(\bar 1-\bar e)\in J(\bar R)$. Let $\bar M\in mspec (\bar R)$. Since $\bar e(\bar 1-\bar e)=\bar 0$,  we have
$mspec( \bar e) \cup mspec(\bar 1-\bar e)=mspec (\bar R)$.  If $\bar M\in  mspec (\bar e)$, we have that $\bar b\in\bar M$ and if $\bar M\in  mspec (\bar 1-\bar e)$ we have that $\bar b(\bar 1-\bar e)\in \bar M$. So, we got that $\bar b(\bar 1-\bar e)\in J(\bar R)$ and $\bar e\in \bar b\bar R$. As a result by [8],  we obtain that $\bar R$ is a semi-regular ring.

Let $\bar R$ be a semi-regular ring. Then for any $\bar b\in \bar R$ there exists such idempotent $\bar e$ that $\bar e\in \bar b\bar R$ and
$\bar b(\bar 1-\bar e)\in J(\bar R)$.

Let $\bar d$ be any non-unit divisor of the element $\bar e$, i.e. $\bar e=\bar d\bar e_1$ for some element $\bar e_1$.
Let $\bar d\bar R+\bar b\bar R=\bar R$. If $\bar d$ is a non-unit element $R$, then there exists maximal ideal $\bar M\in mspec(\bar R)$ such that $\bar d\in \bar M$. Since $\bar b (\bar 1-\bar e)\in J(\bar R)$,  we have $\bar b(\bar 1-\bar e)\in \bar M$. Then $\bar b\in \bar M$ or
$(\bar 1-\bar e)\in \bar M$. If $\bar b\in \bar M$ we have that $\bar 1\in \bar M$ and this is impossible.
If $(\bar 1-\bar e)\in \bar M$ and since $\bar d\in \bar M$,  we obtain that $\bar e=\bar d\bar e_1\in \bar M$. This is impossible because
$\bar e\in \bar b\bar R$. Put $\bar r=\bar 1-\bar e$, $\bar s=\bar e$ we get that $\bar 0=\bar r\bar s$ where
$\bar r\bar R+\bar b\bar R=\bar R$ and
 $\bar s'\bar R+\bar b\bar R\neq\bar R$ for each non-unit divisor $\bar s'$ of $\bar s$.

 Therefore, there exist  such elements $t, u, v \in R$ that $ru+bv=1+at$. Let $aR+rR=\delta R$. Then by Proposition 4,  we have $a=\delta a_0, r=\delta r_0$ and $a_0R+r_0R=R$ for some elements $a_0, r_0\in R$. hence, it follows that $\delta R+bR=R$ and since $\bar r\bar s=\bar 0$,
we see $rs=a\alpha$ for some element $\alpha\in R$. So $\delta r_0s=\delta a_0\alpha$ and  $\delta (r_0s - a_0\alpha)=0$. Since $\delta$ is a nonzero divisor of a  regular element, we have  $r_0s - a_0\alpha=0$. From $a_0R+r_0R=R$ it follows that there exist such elements $k, t\in R$ that
$a_0k+r_0t=1$. This means that $a_0ks+r_0ts=s$, i.e. $a_0\beta =s, \beta \in R$. Thus, $a=\delta a_0$, where  $\delta R+bR=R$ and $sR\subset a_0R$. Let $j$ be such a non-unit divisor of the element $a_0$ that $jR+bR=R$ and $aR\subset jR$.

Hence  $\bar j\bar R+\bar b\bar R=\bar R$,  but we get a contradiction with the condition that $\bar j$ is a non-unit divisor of $\bar s$.
Obtained  contradiction shows that $a=\delta a_0$, where $\delta R+bR=R$ and $a_0'R+bR\neq R$ for any non-unit divisor $a_0'$ of the element $a_0$.

Taking into account that $b$ is an arbitrary element of $R$, we obtain that $a$  is an adequate element of $R$.

The theorem is proved.

As an obvious  consequence of "Test Shores" and Theorem 12, we obtain the following result.

{\tt Theorem 13.} Let $R$ be a semihereditary Bezout ring, in which every regular element  is adequate.
Then $R$  is an elementary divisor ring.

{\it Proof.}  Since a semihereditary Bezout ring is a ring of stable range 2 [13], and for any regular element $a\in R$, the ring $R/aR$ is a semi-regular ring, we obtain that the ring $R/aR$ is a ring of stable range 1 [16]. By [13], $R/aR$ is an elementary divisor ring. And by "Test Shores"
we have that $R$ is an elementary divisor ring.

The theorem is proved.

{\tt Definition 6.} An element $a$ of a ring $R$ is said to be {\it  avoidable} if for any elements $b, c \in R$ such that $aR+bR+cR=R$ there exist such elements $r, s \in R$  that $a=rs$, $rR+bR=R$, $sR+cR=R$, and $rR+sR=R$. A ring $R$ is called {\it  avoidable} if every its nonzero element is  avoidable.

It is easy to see that any adequate element is avoidable.

{\tt Theorem 14.} Let $R$ be a Bezout ring of stable range 2.
A regular element $a\in R$ is  avoidable element iff  $R/aR$ is a clean ring.

{\it Proof.} Let $aR+bR+cR=R, a\neq 0$ and $a=rs$, where  $rR+bR=R$, $sR+cR=R$, and $rR+sR=R$.
Denote $\bar r=r+aR, \bar s=s+aR$. Since $rR+sR=R$,  we have $ru+sv=1$ for some elements $u, v \in R$.
Then $\bar r^2\bar u=\bar r, \bar s^2\bar v =\bar s$. Denote $\bar r\bar u=\bar e$, then $\bar e^2=\bar e$ and  $\bar s\bar v=\bar 1 -\bar e$.

Since $\bar r\bar R+\bar b\bar R=\bar R$, we have $\bar 1 -\bar e \in \bar b\bar  R$.
Similarly, since  $\bar s\bar R+\bar c\bar R=\bar R$, we have $\bar e\in \bar c \bar R$.

We proved that for any $\bar b, \bar c\in \bar R$ there exists an idempotent $\bar e\in \bar R$ such that $\bar e\in \bar b\bar R$ and
$(\bar 1 -\bar e)\in \bar c\bar R$, i,e. $\bar R$ is an exange ring. By [9], $\bar R$ is a clean ring.

Suppose that $\bar R=R/aR$ is a clean ring and $aR+bR+cR=R$.  Denote $\bar b=b+aR, \bar c=c+aR$. Since $\bar R$ is a clean ring, there exists
an idempotent $\bar e\in \bar R$ such that $\bar e\in \bar b\bar R$ and $(\bar 1 -\bar e)\in \bar c\bar R$.
From these conditions we get $e+bp=as$ for some elements $p, s \in R$, and, similarly, $1-e+c\alpha=a\beta$ for some elements $\alpha, \beta \in R$. Since $\bar e=\bar e^2$, we obtain $e(1-e)=at$ for some element $t\in R$.  Let $eR+aR=dR$. By Proposition 4, we have $e=de_0, a=da_0$ and $e_0R+a_0R=R$ for some elements $e_0, a_0\in R$. Then $e+a_0j=1$ for some element $j\in R$. Taking $r=a_0, s=d$, we obtain the decomposition $a=rs$. Since $e+bp=as$, $rR+bR=a_0R+bR=R$, and, similarly, since $1-e+c\alpha=a\beta$,  we have $sR+cR=dR+cR=R$ and obviously
$rR+bR=a_0R+dR=R$.

The theorem is proved.

As in the case of Theorem 13, we obtain the following result.

{\tt Theorem 15.} Let $R$ be a semihereditary Bezout ring in which every regular  element  is avoidable.
Then $R$  is an elementary divisor ring.

Now we will investigate the case when $R$ is a Bezout ring of stable range 2 and $R/aR$ is a semipotent ring.

{\tt Definition 7.}  An element $a$ of a ring $R$ is said to be {\it semipotent} if for any element $b\in R$ such that $b\notin J(aR)$ there exist such non-unit elements $r, s \in R$  that $a=rs$, $rR+bR=R$  and $rR+sR=R$.

Recall  that  a  ring  R  is  a  semipotent  ring,  also  called $I_0$ ring  by
Nicholson [8],  if every principal  ideal not contained  in  $J(R)$  contains  a  nonzero idempotent.  Examples  of  these
rings include: (a)  Exchange and clean rings  (see [9]), (b) Endomorphism rings of injective modules (see [8]),
(c)  Endomorphism  rings  of  regular  modules  in the sense Zelmanowitez (see [18]).

{\tt Theorem 16.} Let $R$ be a Bezout ring of stable range 2.
A regular element $a\in R$ is a semipotent  element iff  $R/aR$ is a semipotent ring.

{\it Proof.} Let  $a=rs$, where  $rR+bR=R$,  and $rR+sR=R$ for each $b\notin J(aR)$.
Since $rR+sR=R$, we have $ru+sv=1$ for some elements $u, v \in R$.
Then $\bar r^2\bar u=\bar r, \bar s^2\bar v =\bar s$. Let  $\bar r\bar u=\bar e$, then $\bar e^2=\bar e$ and  $\bar s\bar v=\bar 1 -\bar e$.

From equality $rR+bR=R$,  it follows that $rx+by=1$ for some elements $x, y\in R$. Then  $\bar 1 -\bar e\in \bar b\bar R$.

Let $\bar R=R/aR$ be a semipotent ring and $\bar b$ be an arbitrary element of $\bar R$ such that $\bar b\notin J(\bar R)$. Then there exists a nontrivial idempotent
$\bar e\in \bar R$ such that $\bar e\in \bar b\bar R$.  Since $\bar e^2=\bar e$,  we have $e(1-e)=as$ for some element $s\in R$. In addition, we have $e-bt=ak$  because  $\bar e\in \bar b\bar R$. Let $aR+eR=dR$.  By Proposition 4, we have $e=de_0, a=da_0$ and $e_0R+a_0R=R$ for some elements $e_0, a_0\in R$. Since the element  $a$ is regular and $e(1-e)=as$,  we obtain $a_0R+eR=R$. Then $a_0p+eq=1$ for some elements $p, q\in R$.
Since $e=ak+bt$,  we have $1=a_0p+eq=a_0(p+dkq)+btq$, i.e. $a_0R+bR=R$. Putting $a_0=r$ and $d=s$, we obtain a desired representation of the element $a$.

The theorem is proved.

Now define the conditions under which the ring $R/aR$ will be a Gelfand ring,  where $R$ is a Bezout ring.

{\tt Definition 8.}  An element $a$ of a ring $R$ is said to be a {\it Gelfand element} if for any elements $b, s\in R$ such that $aR+bR+cR=R$ there exist such  elements $r, s \in R$  that $a=rs$, $rR+bR=R$  and $rR+sR=R$.

Note that according to  [1, 11] we have:

 {\tt Proposition 17.}
Let $R$ be a commutative ring. Then the following are equivalent

1) $R$ is a Gelfand ring.

2) For any elements $a, b\in R$  such that $aR+bR=R$ there exist such elements $c, d \in R$ that $aR+cR=R$, $bR+dR=R$ and $cd=0$.

{\tt Theorem 18.} Let $R$ be a Bezout ring of stable range 2.
A regular element $a\in R$ is a Gelfand element iff  $R/aR$ is a Gelfand ring.

{\it Proof.} Let $a$ is a Gelfand element,  i.e.  $a=rs$, and  $rR+bR=R$, $sR+cR=R$  for any $b, c\in R$ such that  $aR+bR+cR=R$.
Denote $\bar R=R/aR$, $\bar b=b+aR, \bar c=c+aR$  and $\bar r=r+aR, \bar s=s+aR$. It is obvious that $\bar b \bar R+\bar c\bar R=\bar R$. Then
$\bar 0=\bar r\bar c$, and  $\bar r\bar R+\bar b\bar R=\bar R$,  $\bar s\bar R+\bar c\bar R=\bar R$. According to the Proposition 17, we
obtain that $R/aR$ is a Gelfand ring.

If $\bar R=R/aR$ is a  Gelfand ring, we have $\bar 0=\bar r\bar c$, where  $\bar r\bar R+\bar b\bar R=\bar R$,  $\bar s\bar R+\bar c\bar R=\bar R$
for any $b,c\in R$ such that $\bar b\bar R+\bar c\bar R=\bar R$. Hence, we obtain $aR+bR+cR=R$ and $rs\in aR$, where
$\bar b=b+aR, \bar c=c+aR$  and $\bar r=r+aR, \bar s=s+aR$. Then $rs=at$ for some element $t\in R$.

Let $rR+aR=r_1R$, $sR+aR=s_1R$. By Proposition 4, we have $r=r_1r_0, a=r_1a_0$, $s=s_1s_2, a=s_1a_2$ and $r_0R+a_0R=R$,  $s_2R+a_2R=R$.
Since $a$ is a regular element of $R$ and $rs=at$, we obtain $a_0\alpha=s$ for some element $\alpha\in R$. Therefore $a=r_1a_0$, where
$r_1R+bR=R$ and $a_0R+cR=R$.

The theorem is proved.

Note that in proving the necessity of Theorems 12, 14, 16, 18, the condition of the regularity of the element $a$ was not used.

Since the classes of commutative $PM$-rings and Gelfand rings coincide, we obtain the following results.

{\tt Proposition 19.}
Let $R$ be a commutative ring. Then the following statements are equivalent

1) $R$ is a Gelfand ring.

2) For each prime ideal $P\in spec(a)$ there exists a unique maximal ideal $M\in mspec(a)$ such that $P\subset M$.

{\tt Proposition 20.}
The set of all Gelfand elements of a commutative ring is saturated multiplicative closed set.

{\tt Definition 9.} [17] A ring $R$ is said to be a {\it ring of regular range 1} if for any elements $a, b\in R$ such that $aR+bR=R$ there exists such  element $t \in R$  that $a+bt$ is a regular element of $R$.

An obvious example of a ring of regular range 1 is a domain.

Less trivial example of a ring of regular range 1 is a semihereditary ring.

Indeed, let $R$ be a  semihereditary ring and $aR+bR=R$. If $a=0$, we have $bR=R$ and then $a+bt=1$ for some element $t\in R$.
Note that in a semihereditary ring any element can be represented as $er$ where $e$ is an idempotent and $r$ is a regular element of $R$.

Let $a\neq 0$, then for any $t\in R$ we obtain that $a+bt$ is a non-zero element of $R$, i.e. $a+bt=er$, where $e^2=e$  and $r$ is a regular element of $R$. Since $aR+bR=R$, we have $eR+bR=R$. Then $eu+bv=1$  for some elements $u, v \in R$. Then $(1-e)eu+(1-e)bv=1-e$, i.e. $b(1-e)v=1-e$ and $e+bt=1$ for some element $t\in R$. Then $er+btr=r$,  i.e. $R$  is a ring of regular range 1.

{\tt Definition 10.}  A ring $R$ is said to be a {\it ring of avoidable  range 1} if for any elements $a, b\in R$ such that $aR+bR=R$ there exists such  element $t \in R$  that $a+bt$ is an avoidable element of $R$.

{\tt Definition 11.}  A ring $R$ is said to be a {\it ring of Gelfand range 1} if for any elements $a, b\in R$ such that $aR+bR=R$ there exists such  element $t \in R$  that $a+bt$ is a Gelfand element of $R$.

Since any avoidable element is  Gelfand we obtain that a  ring of avoidable  range 1 is a  ring of Gelfand  range 1.
It is obvious that any neat ring  is a  ring of avoidable  range 1 and a $PM^{\ast}$ domain is a  ring of Gelfand  range 1.

{\tt Theorem 21.} Let $R$ be a commutative  Bezout ring of stable range 2 and of Gelfand range 1.
Then $R$  is an elementary divisor ring.

{\it Proof.} Let $A=\begin{pmatrix}a&0\\b&c\end{pmatrix}\ $ and $aR+bR+cR=R$.

For the proof of our statement, according to [3]  it is sufficient to
show that matrix $A$ admits  diagonal reduction.

Let $aR+bR=dR$, i.e. $au+bv=d$ for some elements $u, v\in R$. From the condition $aR+bR+cR=R$ it follows that $bR+dR=R$.
Since $R$ is a ring of Gelfand range 1, we obtain  that $b+(au+bv)t=k$ is a Gelfand element for some element $t$.
Then

$$\begin{pmatrix}1&0\\ut&1\end{pmatrix}\ \begin{pmatrix}a&0\\b&c\end{pmatrix}\
\begin{pmatrix}1&0\\vt&1\end{pmatrix}=\begin{pmatrix}a&0\\k&c\end{pmatrix}\, $$

where, obviously,  $aR+kR+cR=R$ and $r$ is a Gelfand element.

 Then $k=rs$,  where $rR+aR=R$ and $sR+cR=R$. Let $p\in R$ be a such element that $sp+cl=1$ for some element $l\in R$. Hence
$rsp+rcl=r$ and $kp+crl=r$. Denoting $rl=q$, we obtain $(kp+cq)R+aR=R$. Suppose $pR+qR=\delta R$,  i.e.   $p=p_1\delta$, $q=q_1\delta$ and
$\delta =px+qy$, $p_1R+q_1R=R$ for some elements $x, y, p_1, q_1 \in R$. Then from $pR\subset p_1R$ and $pR+cR=R$ $\Rightarrow$ $p_1R+cR=R$, and from $p_1R+q_1R=R$ $\Rightarrow$ $p_1R+(p_1k+q_1c)R=R$.

Since $pk+qc=\delta (p_1k+q_1c)$ and $(pk+qc)R+aR=R$, we obtain $(p_1d+q_1c)R+aR=R$. As well as $p_1R+(p_1d+q_1c)R=R$, finally we have
$p_1aR+(p_1k+q_1c)R=R$. By [3], the matrix $\begin{pmatrix}a&0\\k&c\end{pmatrix}\ $ admits  diagonal reduction. Hence, obviously,
the matrix $\begin{pmatrix}a&0\\b&c\end{pmatrix}\ $ admits  diagonal reduction.

The theorem is proved.

We also get the following result

{\tt Theorem 22.} Let $R$ be a  semihereditary  Bezout ring  and  any regular element of $R$ is Gelfand.
Then $R$  is an elementary divisor ring.

{\it Proof.} Any semihereditary  Bezout ring is a ring of regular range 1 and of stable range 2 [17]. Since  any regular element of $R$ is Gelfand, we have that $R$ is a ring of Gelfand range 1.  By Theorem 21, the theorem is proved.

By (see [15]) we can obtain a new result.

{\tt Theorem 23.}

Let $R$ be a Bezout domain. Then the following statements  are equivalent

1) $R$ is an elementary divisor ring.

2) $R$ is a ring of Gelfand range 1.

{\it Proof.} If $R$ is an elementary divisor ring then by [10] for any elements $a, b, c \in R$ such that $aR+bR=R$ there exists
such element $t\in R$ that $a+bt=uv$, where $uR+cR=R, vR+(1-c)R=R$ and $uR+vR=R$ for some elements $u, v \in R$.
By virtue of the Theorem 14, we have that $a+bt$ is an avoidable element. Since a clean ring is a Gelfand ring, then $R/(a+bt)R$ is
a Gelfand ring, i.e. $R$ is a ring of Gelfand range 1.Theorem 21 completes the proof.

Moreover, we can prove the next result.

{\tt Theorem 24.} Let $R$ be an elementary divisor domain that is not a ring of stable range 1.
Then in $R$ there exists at least one nonunit avoidable element.

{\it Proof.} According to the Theorem 22, we obtain that $R$ is a ring of avoidable  range 1. Since $R$ is not a ring of stable range 1, there exists  at least one nonunit avoidable element.

The Theorem is proved.

\end{document}